\documentclass[12pt]{amsart}
\usepackage[T1]{fontenc}
\usepackage[english]{babel}
\usepackage[english]{babel}
\usepackage{amsmath, amssymb, amsthm, amscd,color,comment}
\usepackage{cancel}
\usepackage{cite}
\usepackage{alltt}
\usepackage[dvipsnames]{xcolor}
\usepackage{array}
\usepackage[small,bf,labelsep=period]{caption}
\usepackage{mathtools}
\usepackage{hyperref}

\usepackage{enumerate}
\newtheorem{theorem}{Theorem}[section]
\newtheorem{definition}[theorem]{Definition}

\newtheorem{remark}[theorem]{Remark}
\newtheorem{lemma}[theorem]{Lemma}
\newtheorem{corollary}[theorem]{Corollary}
\newtheorem{proposition}[theorem]{Proposition}

\DeclarePairedDelimiter\ceil{\lceil}{\rceil}
\DeclarePairedDelimiter\floor{\lfloor}{\rfloor}

\def\la{\lambda}

\def\Q{{\bf Q}}

\def\N{\mathbb N}
\def\R{\mathbb R}
\def\Z{\mathbb Z}
\def\cB{\mathcal B}
\def\cC{\mathcal C}

\def\cH{\mathcal H}

\def\cL{\mathcal L}

\def\cP{\mathcal P}
\def\cM{\mathcal M}

\def\cX{\mathcal X}
\def\cY{\mathcal Y}

\def\fqss{\mathbb F_{q^6}}
\def\fqs{\mathbb F_{q^2}}

\def\fqsn{{\mathbb F}_{q^{2n}}}

\def\fq{\mathbb F_q}

\def\Div{{\rm Div}}

\def\deg{{\rm deg}}

\def\char{{\rm Char}}

\def\char{\mbox{\rm Char}}

\def\negalpha{\text{\boldmath$\alpha$}}

\def\neg1{\text{\boldmath$1$}}
\def\nege{\text{\boldmath$e$}}
\def\negbeta{\text{\boldmath$\beta$}}

\def\neg1{\text{\boldmath$1$}}

\newcommand{\al}{\alpha}
\newcommand{\be}{\beta}

\begin{document}

\title[On generalized Weierstrass Semigroups in arbitrary Kummer extensions]{On generalized Weierstrass Semigroups in arbitrary Kummer extensions of $\fq(x)$}



\thanks{}

\author{Alonso S. Castellanos, Erik A. R. Mendoza, and Guilherme Tizziotti}

\address{Faculdade de Matemática, Universidade Federal de Uberlândia, Campus Santa Mônica, CEP 38400-902, Uberlândia, Brazil}
\email{alonso.castellanos@ufu.br}

\address{Departamento de Matemática, Universidade Federal de Viçosa, Campus Viçosa, CEP 36570-000, Viçosa, Brazil}
\email{erik.mendoza@ufv.br}

\address{Faculdade de Matemática, Universidade Federal de Uberlândia, Campus Santa Mônica, CEP 38400-902, Uberlândia, Brazil}
\email{guilhermect@ufu.br}

\begin{abstract}
In this work, we investigate generalized Weierstrass semigroups in arbitrary Kummer extensions of 
the rational function field $\fq(x)$. We analyze their structure and properties, with a particular emphasis on their maximal elements.
Explicit descriptions of the sets of absolute and relative maximal elements within these semigroups are provided. Additionally, we apply our results to function fields of the maximal curves $\mathcal{X}_{a,b,n,s}$ and $\mathcal{Y}_{n,s}$, which cannot be covered by the Hermitian curve, and the Beelen-Montanucci curve. Our results generalize and unify several earlier contributions in the theory of Weierstrass semigroups, providing new perspectives on the relationship between these semigroups and function fields.
\end{abstract}

\maketitle

\section{Introduction}

Weierstrass semigroups at $n$-tuples of distinct rational places $\Q = (Q_1,\ldots, Q_n)$ in a function field $F$ is a topic of research in algebraic geometry with applications in coding theory. In general, the cases $n=1$ and $n>1$ are analyzed separately. The case $n=1$, known as the Weierstrass semigroup of $F$ at $Q_1$, denoted by $H(Q_1)$, is a classical object in algebraic geometry which is related to many theoretical and applied topics; see e.g. \cite{BMZ2020}, \cite{BLM2021}, \cite{BMV2023}, \cite{C2008}, \cite{GKL1993}, \cite{HLP1998}, \cite{KNT2020}, \cite{E2023} and \cite{S2009}. The study of Weierstrass semigroup of a pair of rational places, $H(Q_1, Q_2)$, was initiated by Arbarello et al., in \cite{ACGH1985}, and the first results about this structure are due to Kim \cite{K1994} and Homma \cite{H1996}. Their work enabled the explicit computation of these semigroups at pairs of rational places for function fields of several curves as Hermitian, Suzuki, Norm-trace, and curves defined by Kummer extensions, which are of interest in coding theory; we refer the reader to \cite{CMQ2016}, \cite{GM2001}, \cite{GM2004}, and \cite{MTT2008}. The general case, at several distinct rational places, was first considered in \cite{CT2005}. Later, in \cite{G2004}, Matthews introduced a notion of minimal generating set of the Weierstrass semigroup $H(\mathbf{Q})$, denoted by $\Gamma(\Q)$, which allows the semigroup to be constructed from a finite subset of elements. Matthews' approach was used to compute the minimal generating set of Weierstrass semigroups at an $n$-tuple of distinct rational places $\Q = (Q_1,\ldots, Q_n)$ in certain function fields, as Hermitian \cite{G2004} and Norm-trace \cite{GD2010}. Since then, distinct techniques have been used to determine the minimal generating set $\Gamma(\Q)$, for $n \geq 2$, on different function fields of curves as the $GK$ curve \cite{CT2018gk}, certain curves of the form $f(y)=g(x)$ \cite{CT2018}, the maximal curves $\mathcal{X}_{a,b,n,s}$ and $\mathcal{Y}_{n,s}$ which cannot be covered by the Hermitian curve \cite{CB2020}, and Kummer extensions, see \cite{BQZ2018} and \cite{YH2017}.

In \cite{D1990}, working with algebraically closed fields, Delgado introduced the concept of \textit{generalized Weierstrass semigroups}, denoted by $\widehat{H}(\Q)$. In \cite{BT2006}, Beelen and Tuta\c{s} studied this concept over finite fields and showed that $H(\Q) = \widehat{H}(\Q) \cap \mathbb{N}_0^n$, where $\mathbb{N}_0$ is the set of nonnegative integers. In the study of generalized Weierstrass semigroups $\widehat{H}(\Q)$, two key objects are particularly important: the sets of absolute maximal elements $\widehat{\Gamma}(\Q)$ and relative maximal elements $\widehat{\Lambda}(\Q)$. Using these sets we can determine the semigroups $\widehat{H}(\Q)$ and, consequently, $H(\Q)$. Studies on these sets and its relations with other objects of the Weierstrass semigroup's theory, as gap sets, pure gap sets, and Riemann-Roch spaces, can be found in \cite{CMT2024}, \cite{CMT2025}, \cite{MTT2019}, \cite{TT2019FF} and \cite{TT2019}. The generalized Weierstrass semigroup for specific function fields was studied in \cite{CMT2025}, \cite{MT2023} and \cite{TT2019}.

In this work, we determine the sets $\widehat{\Gamma}(\Q)$ and $\widehat{\Lambda}(\Q)$ of the generalized Weierstrass semigroup $\widehat{H}(\Q)$, where $\Q$ is an $n$-tuple of distinct pairwise totally ramified places on arbitrary Kummer extensions of the rational function field $\fq(x)$. These results generalize and unify several earlier contributions in the theory of Weierstrass semigroups presented in \cite{CMQ2024}, \cite{CMT2025}, \cite{MT2023}, \cite{TT2019FF}, \cite{TT2019} and \cite{YH2017}. In addition to the generalized Weierstrass semigroup, these sets enable us to identify the minimal generating set $\Gamma(\Q)$ of the Weierstrass semigroup  $H(\Q)$ for an $n$-tuple of rational places not already covered in the literature. 

This paper is organized as follows. In Section \ref{Section 2} we list the definitions and results related to generalized Weierstrass semigroups and Kummer extensions. In Section \ref{Section 3} we provide an arithmetic characterization of the elements in the sets $\widehat{\Gamma}(\Q)$ and $\widehat{\Lambda}(\Q)$, and as a consequence, we obtain an explicit description of these sets. Finally, in Section \ref{Section 4} we determine the sets $\widehat{\Gamma}(\Q)$ and $\widehat{\Lambda}(\Q)$ for the function fields of the maximal curves $\cX_{a, b, n, s}$ and $\cY_{n, s}$, the function field of the curve $y^m=f(x)$, with $f(x)$ separable, and the function field of the Beelen-Montanucci curve.

\section{Preliminaries and Notation} \label{Section 2}

Throughout this article, we let $q$ be the power of a prime $p$ and $\fq$ the finite field with $q$ elements. For $a$ and $b$ integers, we denote by $(a, b)$ the greatest common divisor of $a$ and $b$, by $b \bmod{a} $ the smallest non-negative integer congruent with $b$ modulo $a$, and by $\binom{a}{b}$ the binomial coefficient of $a$ and $b$, with the convention that $\binom{a}{b}=0$ if $a<b$. For $c\in \R$, we denote by $\floor*{c}$, $\ceil*{c}$, and $\{c\}$ the floor, ceiling, and fractional part functions of $c$, respectively. We also denote $\N_0 = \N \cup \{0\}$, where $\N$ is the set of positive integers.

\subsection{Function fields and generalized Weierstrass semigroups} \label{FF GWS}

Let $F/\fq$ be a function field of one variable of genus $g(F)$. We denote by $\mathcal P_{F}$ the set of places in $F$, by $\nu_{P}$ the discrete valuation of $F/\fq$ associated to the place $P\in \cP_{F}$, and by $\Div (F)$ the group of divisors on $F$. For a function $z \in F$, $(z)_{F}, (z)_\infty,$ and $(z)_0$ stands for the principal, pole, and zero divisors of the function $z$ in $F$, respectively. Given a divisor $G\in \Div(F)$ of $F/\fq$, the Riemann-Roch space associated to the divisor $G$ is defined by
$$
\cL(G)=\{z\in F: (z)_{F}+G\geq 0\}\cup \{0\},
$$  
and we denote by $\ell(G)$ its dimension as vector space over $\fq$.

Let $\Q=(Q_1, \dots, Q_n)$ be an $n$-tuple of distinct rational places in $F$. The Weierstrass semigroup $H(\Q)$ and the generalized Weierstrass semigroup $\widehat{H}(\Q)$ of $F$ at $\Q$ are defined, respectively, by the sets
$$
H(\mathbf{Q}) := \left\{(a_{1}, \ldots, a_{n}) \in \mathbb{N}_{0}^ {n} :  (z)_{\infty} = \textstyle\sum_{i=1}^ {n} a_{i}Q_{i} \text{ for some }z\in F \right\}
$$
and
$$
\widehat{H}(\Q):=\{(-\nu_{Q_1}(z),\dots ,-\nu_{Q_n}(z))\in \Z^n : z \in R_\Q\setminus\{0\}\},
$$
where $R_\Q$ denotes the ring of functions on $F$ that are regular outside the places in the set
$\{Q_1,\dots, Q_n\}$. Provided that $q \geq n$, we have that the Weierstrass semigroup of $F$ at $\Q$ can be obtained by the relation $H(\Q)=\widehat{H}(\Q)\cap \N_0^n$, see \cite[Proposition 2]{BT2006}. 

The elements in the finite complement ${G(\mathbf{Q}):=\N^n_0\setminus H(\mathbf{Q})}$ are called \emph{gaps} of $F$ at $\mathbf{Q}$ and can be characterized using Riemann-Roch spaces. In fact, an $n$-tuple $\negalpha=(\al_1, \dots, \al_n)\in \N_0^n$ is an element of $G(\Q)$ if and only if $\ell(D_{\negalpha}(\Q)) = \ell(D_{\negalpha}(\Q) - Q_i)$ for some $1 \leq i \leq n$, where $D_\negalpha(\Q) = \alpha_1 Q_1 + \cdots + \alpha_n Q_n$. A \emph{pure gap} of $F$ at $\Q$ is an $n$-tuple $\negalpha = (\alpha_1, \ldots, \alpha_n) \in G(\mathbf{Q})$ such that $\ell(D_\negalpha(\Q)) = \ell(D_\negalpha(\Q) - Q_i)$ for every $1\leq i \leq n$. The set of pure gaps of $F$ at $\mathbf{Q}$ will be denoted by $G_0(\mathbf{Q})$.

The elements of the generalized semigroup $\widehat{H}(\Q)$ also can be characterized using Riemann-Roch spaces as follows. 

\begin{proposition}\cite[Proposition 2.2]{MTT2019} \label{prop_hat_H(Q)}
Let $\negalpha \in \mathbb{Z}^n$ and assume that $q \geq n$. Then
$$\negalpha \in \widehat{H} (\mathbf{Q})\quad \text{if and only if}\quad \cL(D_{\negalpha}(\Q)) \neq \cL(D_{\negalpha}(\Q) - Q_i) \text{ for every } 1\leq i\leq n.$$
\end{proposition}

A fundamental concept in the study of generalized Weierstrass semigroups $\widehat{H}(\Q)$, originally introduced by Delgado \cite{D1990}, is notion of maximal elements, which play a crucial role in describing the elements of $\widehat{H}(\Q), G(\Q)$ and $G_0(\Q)$, see \cite{CMT2024}, \cite{CMT2025}, \cite{MTT2019}, and \cite{TT2019}.

To introduce the concept of maximal elements we need to define the sets below. Set $I:=\{1,\ldots,n\}$. For $i\in I$, a nonempty subset $J\subsetneq I$, and $\negalpha=(\alpha_1,\ldots,\alpha_n)\in \Z^n$, we shall denote
\begin{itemize}
\item [$\bullet$] $\overline{\nabla}_J (\negalpha):=\{\negbeta\in \Z^n : \be_j=\al_j \text{ for }j\in J \text{ and } \be_i < \al_i \text{ for }i\notin J\}$,
\item [$\bullet$] $\nabla_J(\negalpha):=\overline{\nabla}_J(\negalpha)\cap \widehat{H}(\Q)$,
\item [$\bullet$] $\overline{\nabla}(\negalpha):=\cup_{i=1}^{n}\overline{\nabla}_i(\negalpha)$, where $\overline{\nabla}_i(\negalpha):= \overline{\nabla}_{\{i\}}(\negalpha)$, and
\item [$\bullet$] $\nabla(\negalpha):=\overline{\nabla}(\negalpha)\cap \widehat{H}(\Q)$.
\end{itemize}

\begin{definition} \label{defi maximals}
An element $\negalpha\in \widehat{H}(\Q)$ is called maximal if $\nabla(\negalpha)=\emptyset$. If moreover $\nabla_J(\negalpha)=\emptyset$ for every $J\subsetneq I$ with $|J| \geq 2$, we say that $\negalpha$ is absolute maximal. If $\negalpha$ is maximal and $\nabla_J(\negalpha)\neq\emptyset$ for every $J\subsetneq I$ with $|J| \geq 2$, we say that $\negalpha\in \widehat{H}(\Q)$ is relative maximal. The sets of absolute and relative maximal elements in $\widehat{H}(\Q)$ will be denoted, respectively, by $\widehat{\Gamma}(\Q)$ and $\widehat{\Lambda}(\Q)$.
\end{definition}

Observe that the notions of absolute and relative maximality coincide when $n=2$. Next, we present results that characterize the absolute and relative maximal elements of $\widehat{H}(\Q)$. To do this, we introduce the concept of discrepancy, as defined by Duursma and Park in \cite{DP2012}.



\begin{definition} \label{discrepancy}
Let $P_1$ and $P_2$ be distinct rational places on $F$. A divisor $A\in \Div(F)$ is called a discrepancy with respect to $P_1$ and $P_2$ if
$$
\cL(A)\neq \cL(A-P_1)\quad \text{and}\quad\cL(A-P_2)=\cL(A-P_1-P_2).
$$ 
\end{definition}



The following results establishes equivalences using the notion of discrepancy.

\begin{proposition}\cite[Proposition 3]{TT2019FF} \label{prop_discrepancia_Gamma}
Let $\negalpha \in \mathbb{Z}^n$ and assume that $q\geq n$. The following statements are equivalent:
\begin{enumerate}[\rm (i)]
\item $\negalpha \in \widehat{\Gamma}(\mathbf{Q})$;
\item $D_{\negalpha}(\mathbf{Q})$ is a discrepancy with respect to any pair of distinct places in $\{Q_1,\ldots,Q_{n}\}$.
\end{enumerate}
\end{proposition}

\begin{proposition}\cite[Proposition 2.8]{TT2019} \label{prop_discrepamcia_Lambda}
Let $\negalpha\in \mathbb{Z}^n$ and assume that $q \geq n$. The following statements are equivalent:
\begin{enumerate}[\rm (i)]
\item $\negalpha\in \widehat{\Lambda}(\Q)$;
\item $D_{\negalpha-\textbf{1}}(\Q)+Q_i+Q_j$ is a discrepancy with respect to $Q_i$ and $Q_j$ for every $i,j\in I$ with $i\neq j$, where $\neg1$ is the $n$-tuple whose coordinates are all $1$.
\end{enumerate}
\end{proposition}

\subsection{Kummer extensions} \label{Kummer}
Let $m\geq 2$ be a integer such that $\char(\fq) \nmid m$ and $f(x)\in \fq[x]$ be a polynomial with $\deg(f)\geq 2$ such that $f(x)$ is not a $d$-th power (where $d$ divides $m$) of an element in $\fq(x)$. Assume that all roots of $f(x)$ are in $\fq$. Consider the algebraic curve defined by the affine equation
\begin{equation}\label{equation}
\cX: \quad y^m=f(x)
\end{equation}
with function field $\fq(\cX)$. Then $\fq(\cX)/\fq(x)$ is a Kummer extension. Let
\begin{itemize}
\item $P_1, P_2, \dots, P_r\in \cP_{\fq(x)}$ be the zeros and poles of $f(x)$,
\item $\la_{k}:=\nu_{P_k}(f(x))$ be the multiplicity of $P_k$ in $f(x)$, 
\item $P_1, P_2, \dots, P_n\in \cP_{\fq(x)}$ (where $n\leq r$) be distinct pairwise totally ramified places in the extension $\fq(\cX)/\fq(x)$,
\item $Q_k\in \cP_{\fq(\cX)}$ be the unique extension in $\fq(\cX)$ over $P_k$ for $k=1, \dots, n$, and
\item $\Q=(Q_1, Q_2, \dots, Q_n)$.
\end{itemize} 

In \cite[Theorem 3.3]{BMMQ2018}, using a result by Maharaj \cite{M2004}, the authors provided an arithmetic criterion to determine when an element $\negalpha\in \N_0^n$ is a gap of $\fq(\cX)$ at $\Q$. A simplified version of this result is presented below.
\begin{theorem}\cite[Theorem 3.3]{BMMQ2018}\label{teo_gap}
Let $\negalpha=(\al_1, \dots, \al_n)\in \N_0^n$ and $i\in I$. Then $\cL(D_{\negalpha}(\Q))=\cL(D_{\negalpha}(\Q)-Q_i)$ if and only if for every $t\in \{0, \dots, m-1\}$ exactly one of the two following conditions is satisfied:
\begin{enumerate}[\rm (i)]
\item $\displaystyle\sum_{k=1}^{n}\floor*{\frac{\al_k+t\la_k}{m}}+\sum_{k=n+1}^{r}\floor*{\frac{t\la_k}{m}}<0,$
\vspace{2mm}
\item $\displaystyle\floor*{\frac{\al_i+t\la_i}{m}}=\floor*{\frac{\al_i-1+t\la_i}{m}}.
$
\end{enumerate}
\end{theorem}
It is not difficult to verify that, using the same proof presented in \cite[Theorem 3.3]{BMMQ2018}, this arithmetic criterion remains valid for any element $\negalpha \in \Z^n$. On the other hand, a very useful arithmetic result that will be used in this work is the following.
\begin{lemma}\cite[Lemma 3.3]{ABQ2019}\label{lemma_floor}
Let $\al, \la\in \Z$ and $m\in \N$ be such that $(\la, m)=1$. Then there exists a unique integer $t\in \{0, \dots, m-1\}$ such that $\al+t\la\equiv 0 \bmod{m}$ and, for $i\in \Z$, we have
$$
\floor*{\frac{\al+i\la}{m}}=\left\{
\begin{array}{ll}
\displaystyle\floor*{\frac{\al-1+i\la}{m}}+1=\ceil*{\frac{\al}{m}}+\floor*{\frac{i\la}{m}}, & \text{if } i\equiv t \bmod{m},\\
\\
\displaystyle\floor*{\frac{\al-1+i\la}{m}}, & \text{if } i\not\equiv t \bmod{m}.
\end{array}\right.
$$  
\end{lemma}
Using this lemma and considering that the result given in Theorem \ref{teo_gap} is valid for any element $\negalpha\in \Z^n$, we can present a revised version of this result as follows.
\begin{theorem}\label{teo_equiv_gap}
Let $\negalpha=(\al_1, \dots, \al_n) \in \Z^n$ and $i\in I$. Then $\cL(D_{\negalpha}(\Q))=\cL(D_{\negalpha}(\Q)-Q_i)$ if and only if
$$
\sum_{k=1}^{n}\floor*{\frac{\al_k+t_i\la_k}{m}}+\sum_{k=n+1}^{r}\floor*{\frac{t_i\la_k}{m}}<0,
$$
where $t_i\in\{0, \dots, m-1\}$ is the unique integer such that $\al_i+t_i\la_i\equiv 0 \bmod{m}$.
\end{theorem}

In \cite{CMS2025}, the authors provided an explicit description of the gap set $G(P)$ at any totally ramified place $P\in \cP_{\fq(\cX)}$ in the extension $\fq(\cX)/\fq(x)$ through the following functions. For each $1 \leq i \leq m-1$ and $1 \leq k \leq n$, set
\begin{equation} \label{t beta}
t_k(i):= (i\la_k)\bmod m\quad \text{and}\quad \be(i):=\sum_{k=1}^{r}\ceil*{\frac{i\la_k}{m}}-1.
\end{equation}

\begin{remark} \label{t diferente de zero}
Note that, for each $1 \leq k \leq n$, since $P_k$ is a totally ramified place, we have that $(\la_k , m)=1$, and thus $t_k(i) \in \{1, \ldots, m-1\}$ for every $1 \leq i \leq m-1$.
\end{remark}




\section{On maximal elements in generalized Weierstrass semigroups} \label{Section 3}

Let $\cX$ be the algebraic curve defined in (\ref{equation}), $\fq(\cX)/\fq(x)$ be the Kummer extension, and $\Q=(Q_1, Q_2, \dots, Q_n)$ be the $n$-tuple of places in $ \cP_{\fq(\cX)}$ as in Section \ref{Kummer}. In this section, using the functions $t_k(i)$ and $\beta(i)$ given in (\ref{t beta}), we describe $\widehat{\Gamma}(\Q)$ and $\widehat{\Lambda}(\Q)$, the sets of absolute and relative maximals elements in $\widehat{H}(\Q)$, respectively. 

In Section \ref{FF GWS}, we observed that the concept of discrepancy is a very useful tool for determining the sets $\widehat{\Gamma}(\Q)$ and $\widehat{\Lambda}(\Q)$. We know that for a divisor to be a discrepancy with respect to two places, two conditions involving Riemann-Roch spaces must be satisfied, which in general is somewhat difficult to verify. The next results provide an  arithmetic criterion for a certain divisor to be a discrepancy.

\begin{lemma}\label{lemma_discrepancy_Kummer}
Assume that $2\leq n \leq q$ and let $\negalpha = (\alpha_1, \ldots, \alpha_n) \in \Z^n$. Then $D_\negalpha(\Q)$ is a discrepancy with respect to any pair of distinct places in $\{Q_1, Q_2, \dots, Q_n\}$ if and only if there exists a unique $t\in \{0, \dots, m-1\}$ such that $\al_k+t\la_k\equiv 0 \bmod{m}$ for every $1\leq k \leq n$ and
$$
\sum_{k=1}^{n}\ceil*{\frac{\al_k}{m}}+\sum_{k=1}^{r}\floor*{\frac{t\la_k}{m}}=0.
$$
\end{lemma}
\begin{proof}
Let $2\leq n \leq q$ and $\negalpha \in \Z^n$. We begin by noting that, from Definition \ref{discrepancy} and Proposition \ref{prop_hat_H(Q)}, $D_\negalpha(\Q)$ is a discrepancy with respect to any pair of distinct places in $\{Q_1, Q_2, \dots, Q_n\}$ if and only if $\negalpha \in \widehat{H}(\Q)$ and $\cL(D_{\negalpha-\nege_j}(\Q))=\cL(D_{\negalpha-\nege_j}(\Q)-Q_i)$ for every $i, j \in I$ with $i \neq j$, where $\nege_j$ is the $n$-tuple whose the $j$-th coordinate is $1$ and the others are $0$. Furthermore, from Proposition \ref{prop_hat_H(Q)} and Theorem \ref{teo_equiv_gap}, we have that

$\bullet$ $\negalpha \in \widehat{H}(\Q)$ if and only if
\begin{equation}\label{equation_1}
\sum_{k=1}^{n}\floor*{\frac{\al_k+t_i\la_k}{m}}+\sum_{k=n+1}^{r}\floor*{\frac{t_i\la_k}{m}}\geq 0\quad \text{for every } i\in I,
\end{equation} 
where $t_i\in \{0, \dots, m-1\}$ is the unique integer such that $\al_i+t_i\la_i\equiv 0\bmod{m}$, and 

$\bullet$ $\cL(D_{\negalpha-\nege_j}(\Q))=\cL(D_{\negalpha-\nege_j}(\Q)-Q_i)$ for every $i, j \in I$ with $i \neq j$ if and only if
\begin{equation}\label{equation_2}
\sum_{\substack{k=1\\k\neq j}}^{n}\floor*{\frac{\al_k+t_i\la_k}{m}}+\floor*{\frac{\al_j-1+t_i\la_j}{m}}+\sum_{k=n+1}^{r}\floor*{\frac{t_i\la_k}{m}}< 0\quad  \text{for every } i, j\in I \text{ with } i\neq j. 
\end{equation}

Now, suppose that $D_\negalpha(\Q)$ is a discrepancy with respect to any pair of distinct places in $\{Q_1, Q_2, \dots, Q_n\}$. So, Equations (\ref{equation_1}) and  (\ref{equation_2}) are satisfied. This implies that
\begin{equation}\label{equation_3}
\sum_{k=1}^{n}\floor*{\frac{\al_k+t_i\la_k}{m}}+\sum_{k=n+1}^{r}\floor*{\frac{t_i\la_k}{m}}=0 \quad \text{for every }i\in I 
\end{equation}
and
\begin{equation}\label{equation_4}
\floor*{\frac{\al_j+t_i\la_j}{m}}=\floor*{\frac{\al_j-1+t_i\la_j}{m}}+1\quad  \text{for every } i, j\in I \text{ with } i\neq j.
\end{equation}
From Lemma \ref{lemma_floor} and Equation (\ref{equation_4}) we obtain $t:=t_1=t_2=\dots=t_n$. Thus, $\al_k+t\la_k\equiv 0 \bmod{m}$ for every $1\leq k \leq n$. Furthermore, from Lemma \ref{lemma_floor} and Equation (\ref{equation_3}) we conclude that
\begin{align*}
\sum_{k=1}^{n}\ceil*{\frac{\al_k}{m}}+\sum_{k=1}^{r}\floor*{\frac{t\la_k}{m}}&=\sum_{k=1}^{n}\left(\floor*{\frac{\al_k+t\la_k}{m}}-\floor*{\frac{t\la_k}{m}}\right)+\sum_{k=1}^{r}\floor*{\frac{t\la_k}{m}}\\
&=\sum_{k=1}^{n}\floor*{\frac{\al_k+t\la_k}{m}}+\sum_{k=n+1}^{r}\floor*{\frac{t\la_k}{m}}\\
&=0.
\end{align*}

Conversely, suppose that there exists a unique $t\in \{0, \dots, m-1\}$ such that $\al_k+t\la_k\equiv 0 \bmod{m}$ for every $1\leq k \leq n$ and
$$
\sum_{k=1}^{n}\ceil*{\frac{\al_k}{m}}+\sum_{k=1}^{r}\floor*{\frac{t\la_k}{m}}=0.
$$ 
From Lemma \ref{lemma_floor} we get 
$$
\sum_{k=1}^{n}\floor*{\frac{\al_k+t\la_k}{m}}+\sum_{k=n+1}^{r}\floor*{\frac{t\la_k}{m}}=\sum_{k=1}^{n}\ceil*{\frac{\al_k}{m}}+\sum_{k=1}^{r}\floor*{\frac{t\la_k}{m}}=0
$$ 
and, for every $j\in I$,
\begin{align*}
&\sum_{\substack{k=1\\k\neq j}}^{n}\floor*{\frac{\al_k+t\la_k}{m}}+\floor*{\frac{\al_j-1+t\la_j}{m}}+\sum_{k=n+1}^{r}\floor*{\frac{t\la_k}{m}}\\
&=\sum_{\substack{k=1\\k\neq j}}^{n}\floor*{\frac{\al_k+t\la_k}{m}}+\floor*{\frac{\al_j+t\la_j}{m}}-1+\sum_{k=n+1}^{r}\floor*{\frac{t\la_k}{m}}\\
&=\sum_{k=1}^{n}\floor*{\frac{\al_k+t\la_k}{m}}-1+\sum_{k=n+1}^{r}\floor*{\frac{t\la_k}{m}}\\
&=-1.
\end{align*} 
This implies that $\negalpha$ satisfies Equations (\ref{equation_1}) and (\ref{equation_2}). Therefore, $D_\negalpha(\Q)$ is a discrepancy with respect to any pair of distinct places in $\{Q_1, Q_2, \dots, Q_n\}$. 
\end{proof}

\begin{lemma}\label{lemma2_discrepancy_Kummer}
Assume that $2 \leq n \leq q$ and let $\negalpha=(\al_1, \dots, \al_n) \in \Z^n$. Then $D_{\negalpha-\textbf{1}}(\Q)+Q_i+Q_j$ is a discrepancy with respect to $Q_i$ and $Q_j$ for every $i,j\in I$ with $i\neq j$ if and only if there exists a unique $t\in \{0, \dots, m-1\}$ such that $\al_k+t\la_k\equiv 0 \bmod{m}$ for every $1\leq k \leq n$ and
$$
\sum_{k=1}^{n}\ceil*{\frac{\al_k}{m}}+\sum_{k=1}^{r}\floor*{\frac{t\la_k}{m}}=n-2.
$$
\end{lemma}
\begin{proof}
From Definition \ref{discrepancy} the divisor $D_{\negalpha-\neg1}(\Q)+Q_i+Q_j$ is a discrepancy with respect to $Q_i$ and $Q_j$ for every $i, j\in I$ with $i\neq j$ if and only if $\cL(D_{\negalpha-\neg1}(\Q)+Q_i+Q_j)\neq \cL(D_{\negalpha-\neg1}(\Q)+Q_j)$ for every $i, j \in I$ with $i \neq j$ and $\cL(D_{\negalpha-\neg1}(\Q)+Q_i)=\cL(D_{\negalpha-\neg1}(\Q))$ for every $i\in I$. Furthermore, from Theorem \ref{teo_equiv_gap}, we have that

$\bullet$ $\cL(D_{\negalpha-\neg1}(\Q)+Q_i+Q_j)\neq \cL(D_{\negalpha-\neg1}(\Q)+Q_j)$ for every $i, j \in I$ with $i \neq j$ if and only if
\begin{equation}\label{equation_5}
\sum_{\substack{k=1\\k\neq i, j}}^{n}\floor*{\frac{\al_k-1+t_i\la_k}{m}}+\floor*{\frac{\al_i+t_i\la_i}{m}}+\floor*{\frac{\al_j+t_i\la_j}{m}}+\sum_{k=n+1}^{r}\floor*{\frac{t_i\la_k}{m}}\geq 0
\end{equation}
for every $i, j\in I$ with $i\neq j$, where $t_i\in \{0, \dots, m-1\}$ is the unique integer such that $\al_i+t_i\la_i\equiv 0\bmod{m}$, and

$\bullet$ $\cL(D_{\negalpha-\neg1}(\Q)+Q_i)=\cL(D_{\negalpha-\neg1}(\Q))$ for every $i\in I$ if and only if
\begin{equation}\label{equation_6}
\sum_{\substack{k=1\\k\neq i}}^{n}\floor*{\frac{\al_k-1+t_i\la_k}{m}}+\floor*{\frac{\al_i+t_i\la_i}{m}}+\sum_{k=n+1}^{r}\floor*{\frac{t_i\la_k}{m}}< 0\quad  \text{for every } i\in I. 
\end{equation}

Suppose that $D_{\negalpha-\textbf{1}}(\Q)+Q_i+Q_j$ is a discrepancy with respect to $Q_i$ and $Q_j$ for every $i,j\in I$ with $i\neq j$. So, Equations (\ref{equation_5}) and (\ref{equation_6}) are satisfied. This implies 
\begin{equation}\label{equation_7}
\sum_{\substack{k=1\\k\neq i, j}}^{n}\floor*{\frac{\al_k-1+t_i\la_k}{m}}+\floor*{\frac{\al_i+t_i\la_i}{m}}+\floor*{\frac{\al_j+t_i\la_j}{m}}+\sum_{k=n+1}^{r}\floor*{\frac{t_i\la_k}{m}}=0
\end{equation}
and
\begin{equation}\label{equation_8}
\floor*{\frac{\al_j+t_i\la_j}{m}}=\floor*{\frac{\al_j-1+t_i\la_j}{m}}+1
\end{equation}
for every $i, j\in I$ with $i\neq j$. From Lemma \ref{lemma_floor} and Equation (\ref{equation_8}) we obtain $t:=t_1=t_2=\dots=t_n$ and thus $\al_k+t\la_k\equiv 0 \bmod{m}$ for every $1\leq k \leq n$. Moreover, from Lemma \ref{lemma_floor} and Equation (\ref{equation_7}), we can conclude that
\begin{align*}
\sum_{k=1}^{n}\ceil*{\frac{\al_k}{m}}+\sum_{k=1}^{r}\floor*{\frac{t\la_k}{m}}&=\sum_{k=1}^{n}\left(\floor*{\frac{\al_k+t\la_k}{m}}-\floor*{\frac{t\la_k}{m}}\right)+\sum_{k=1}^{r}\floor*{\frac{t\la_k}{m}}\\
&=\sum_{k=1}^{n}\floor*{\frac{\al_k+t\la_k}{m}}+\sum_{k=n+1}^{r}\floor*{\frac{t\la_k}{m}}\\
&=\sum_{\substack{k=1\\k\neq i, j}}^{n}\floor*{\frac{\al_k+t\la_k}{m}}+\floor*{\frac{\al_i+t\la_i}{m}}+\floor*{\frac{\al_j+t\la_j}{m}}+\sum_{k=n+1}^{r}\floor*{\frac{t\la_k}{m}}\\
&=\sum_{\substack{k=1\\k\neq i, j}}^{n}\left(\floor*{\frac{\al_k-1+t\la_k}{m}}+1\right)+\floor*{\frac{\al_i+t\la_i}{m}}+\floor*{\frac{\al_j+t\la_j}{m}}+\sum_{k=n+1}^{r}\floor*{\frac{t\la_k}{m}}\\
&=n-2+\sum_{\substack{k=1\\k\neq i, j}}^{n}\floor*{\frac{\al_k-1+t\la_k}{m}}+\floor*{\frac{\al_i+t\la_i}{m}}+\floor*{\frac{\al_j+t\la_j}{m}}+\sum_{k=n+1}^{r}\floor*{\frac{t\la_k}{m}}\\
&=n-2.
\end{align*}

Conversely, suppose that there exists a unique $t\in \{0, \dots, m-1\}$ such that $\al_k+t\la_k\equiv 0 \bmod{m}$ for every $1\leq k \leq n$ and
$$
\sum_{k=1}^{n}\ceil*{\frac{\al_k}{m}}+\sum_{k=1}^{r}\floor*{\frac{t\la_k}{m}}=n-2.
$$ 
Then, from Lemma \ref{lemma_floor}, for every $i, j\in I$ with $i\neq j$  we have that 
\begin{align*}
&\sum_{\substack{k=1\\k\neq i, j}}^{n}\floor*{\frac{\al_k-1+t\la_k}{m}}+\floor*{\frac{\al_i+t\la_i}{m}}+\floor*{\frac{\al_j+t\la_j}{m}}+\sum_{k=n+1}^{r}\floor*{\frac{t\la_k}{m}}\\
&=\sum_{\substack{k=1\\k\neq i, j}}^{n}\left(\ceil*{\frac{\al_k}{m}}+\floor*{\frac{t\la_k}{m}}-1\right)+\ceil*{\frac{\al_i}{m}}+\floor*{\frac{t\la_i}{m}}+\ceil*{\frac{\al_j}{m}}+\floor*{\frac{t\la_j}{m}}+\sum_{k=n+1}^{r}\floor*{\frac{t\la_k}{m}}\\
&=-n+2+\sum_{k=1}^{n}\ceil*{\frac{\al_k}{m}}+\sum_{k=1}^{r}\floor*{\frac{t\la_k}{m}}\\
&=0
\end{align*}
and, for every $i\in I$,
\begin{align*}
&\sum_{\substack{k=1\\k\neq i}}^{n}\floor*{\frac{\al_k-1+t\la_k}{m}}+\floor*{\frac{\al_i+t\la_i}{m}}+\sum_{k=n+1}^{r}\floor*{\frac{t\la_k}{m}}\\
&=\sum_{\substack{k=1\\k\neq i}}^{n}\left(\ceil*{\frac{\al_k}{m}}+\floor*{\frac{t\la_k}{m}}-1\right)+\ceil*{\frac{\al_i}{m}}+\floor*{\frac{t\la_i}{m}}+\sum_{k=n+1}^{r}\floor*{\frac{t\la_k}{m}}\\
&=-n+1+\sum_{k=1}^{n}\ceil*{\frac{\al_k}{m}}+\sum_{k=1}^{r}\floor*{\frac{t\la_k}{m}}\\
&=-1.
\end{align*} 
This implies that $\negalpha$ satisfy Equations (\ref{equation_5}) and (\ref{equation_6}). Therefore, $D_{\negalpha-\textbf{1}}(\Q)+Q_i+Q_j$ is a discrepancy with respect to $Q_i$ and $Q_j$ for every $i,j\in I$ with $i\neq j$.
\end{proof}

For abbreviation, in the remainder of this section, $\Upsilon$ will denote either $\Gamma$ or $\Lambda$. Using this notation and Propositions \ref{prop_discrepancia_Gamma} and \ref{prop_discrepamcia_Lambda}, we can summarize the lemmas above as follows.

\begin{lemma}\label{lemma_criterion_Upsilon}
Assume that $2 \leq n \leq q$ and let $\negalpha=(\al_1, \dots, \al_n)\in \Z^n$. Then $\negalpha\in \widehat{\Upsilon}(\Q)$ if and only if there exists a unique $t\in \{0, \dots, m-1\}$ such that $\al_k+t\la_k\equiv 0 \bmod{m}$ for every $1\leq k \leq n$ and 
$$
\sum_{k=1}^{n}\ceil*{\frac{\al_k}{m}}+\sum_{k=1}^{r}\floor*{\frac{t\la_k}{m}}=\rho,\quad \text{where}\quad \rho=\left\{
\begin{array}{ll}
	0, & \text{if }\, \Upsilon=\Gamma,\\
	n-2, & \text{if }\, \Upsilon=\Lambda.
\end{array}\right.
$$  
\end{lemma}

We can now present our main result, in which we provide an explicit description of the set $\widehat{\Upsilon}(\Q)$.

\begin{theorem}\label{teo_widehat_Upsilon}
Let $2\leq n \leq q$ and $\Q=(Q_{1}, Q_{2}, \dots, Q_{n})$ be as above. Then
\begin{align*}
\widehat{\Upsilon}(\Q)=&\Bigg\{(mj_1+t_{1}(i), \dots,  mj_n+t_{n}(i)):
\begin{array}{l}
	1\leq i \leq m-1,\,\, j_1, \dots, j_n\in \Z, \\
	j_1+\cdots +j_n=\beta(i)+1-n+\rho
\end{array}
\Bigg\}\\
& \bigcup \Bigg\{(mj_1, mj_2, \dots, mj_n):\, j_1, j_2, \dots, j_n \in \Z, \, \,   j_1+j_2+\cdots+j_n=\rho \Bigg\},
\end{align*}
where 
$$\quad \rho=\left\{
\begin{array}{ll}
0, & \text{if }\, \Upsilon=\Gamma,\\
n-2, & \text{if }\, \Upsilon=\Lambda.
\end{array}\right.
$$  
\end{theorem}
\begin{proof}
Define the set 
\begin{align*}
	\widehat{\Upsilon}=&\Bigg\{(mj_1+t_{1}(i), \dots,  mj_n+t_{n}(i)):
	\begin{array}{l}
		1\leq i \leq m-1,\,\, j_1, \dots, j_n\in \Z, \\
		j_1+\cdots +j_n=\beta(i)+1-n+\rho
	\end{array}
	\Bigg\}\\
	& \bigcup \Bigg\{(mj_1, mj_2, \dots, mj_n):\, j_1, j_2, \dots, j_n \in \Z, \, \,   j_1+j_2+\cdots+j_n=\rho \Bigg\}.
\end{align*}
Let $\negalpha\in \widehat{\Upsilon}$. We will prove that $\negalpha\in \widehat{\Upsilon}(\Q)$. To do this, we will analyze two cases. 

$\bullet$ Case $\negalpha=(mj_1+t_1(i), \dots, mj_n+t_n(i))$. For this case we have that there exists a unique integer $t:=m-i\in \{0, \dots, m-1\}$ such that $mj_k+t_k(i)+t\la_k\equiv 0 \bmod{m}$ for every $1\leq k \leq n$. Furthermore, from definition of $\la_k$ given in Section \ref{Kummer}, we have $\sum_{k=1}^{r}\la_k=0$, and from Remark \ref{t diferente de zero} we have $1 \leq t_{k}(i) \leq m-1$. So, we obtain
\begin{align*}
\sum_{k=1}^{n}\ceil*{\frac{mj_k+t_k(i)}{m}}+\sum_{k=1}^{r}\floor*{\frac{t\la_k}{m}}&=
\sum_{k=1}^{n}\left(j_k+\ceil*{\frac{t_k(i)}{m}}\right)+\sum_{k=1}^{r}\floor*{\frac{(m-i)\la_k}{m}}\\
&=j_1+j_2+\cdots+j_n+n+\sum_{k=1}^{r}\la_k-\sum_{k=1}^{r}\ceil*{\frac{i\la_k}{m}}\\
&=\be(i)+1+\rho-\sum_{k=1}^{r}\ceil*{\frac{i\la_k}{m}}\\
&=\rho.
\end{align*}

$\bullet$ Case $\negalpha=(mj_1, \dots, mj_n)$. As in the previous case, we can conclude that there exists a unique integer $t:=0\in \{0, \dots, m-1\}$ such that $mj_k+t\la_k\equiv 0 \bmod{m}$ for every $1\leq k \leq n$. Furthermore, 
$$
\sum_{k=1}^{n}\ceil*{\frac{mj_k}{m}}+\sum_{k=1}^{r}\floor*{\frac{t\la_k}{m}}=j_1+j_2+\cdots+j_n=\rho.
$$
Thus, by Lemma \ref{lemma_criterion_Upsilon}, we get $\negalpha\in \widehat{\Upsilon}(\Q)$ and consequently $\widehat{\Upsilon}\subseteq \widehat{\Upsilon}(\Q)$. 

Conversely, let $\negalpha=(\al_1, \dots, \al_n) \in \widehat{\Upsilon}(\Q)$. From Lemma \ref{lemma_criterion_Upsilon}, there exists a unique $i\in \{0, \dots, m-1\}$ such that $\al_k+i\la_k\equiv 0 \bmod{m}$ for every $1\leq k\leq n$, and
$$
\sum_{k=1}^{n}\ceil*{\frac{\al_k}{m}}+\sum_{k=1}^{r}\floor*{\frac{i\la_k}{m}}=\rho.
$$
For each $1\leq k \leq n$, we can write $\al_k=mj_k+i_k$, where $j_k\in \Z$ and $0\leq i_k \leq m-1$. 

$\bullet$ Case $i=0$. In this case we have that $\al_k\equiv 0 \bmod{m}$ for every $1\leq k\leq n$, and thus $i_1=i_2=\dots=i_n=0$. So, 
$$
\rho=\sum_{k=1}^{n}\ceil*{\frac{\al_k}{m}}+\sum_{k=1}^{r}\floor*{\frac{i\la_k}{m}}=\sum_{k=1}^{n}\ceil*{\frac{mj_k}{m}}=j_1+j_2+\cdots+j_n.
$$
Therefore, $\negalpha=(mj_1, mj_2, \dots, mj_n) \in \widehat{\Upsilon}$.

$\bullet$ Case $i\neq 0$. From the relations $\al_k+i\la_k\equiv 0 \bmod{m}$ and $\al_k=mj_k+i_k$ we obtain that $i_k \equiv -i \la_k \bmod{m}$ and therefore $i_k \equiv (m-i)\la_k \bmod{m}$. Now, since $i\neq 0$ we have $1\leq i_k \leq m-1$. Thus, from Equation (\ref{t beta}) and Remark \ref{t diferente de zero}, we conclude that $i_k=t_k(m-i)$ for every $1\leq k \leq n$. So,
\begin{align*}
\rho&=\sum_{k=1}^{n}\ceil*{\frac{\al_k}{m}}+\sum_{k=1}^{r}\floor*{\frac{i\la_k}{m}}\\
&=\sum_{k=1}^{n}\ceil*{\frac{mj_k+t_k(m-i)}{m}}-\sum_{k=1}^{r}\ceil*{\frac{(m-i)\la_k}{m}}\\
&=\sum_{k=1}^{n}\left(j_k+\ceil*{\frac{t_k(m-i)}{m}}\right)-\sum_{k=1}^{r}\ceil*{\frac{(m-i)\la_k}{m}}\\
&=j_1+j_2+\cdots+j_n+n-\be(m-i)-1.
\end{align*}
This implies that $\negalpha=(mj_1+t_1(m-i), \dots, mj_n+t_n(m-i)) \in \widehat{\Upsilon}$, and consequently $\widehat{\Upsilon}(\Q)\subseteq \widehat{\Upsilon}$. 
\end{proof}

The set $\Gamma(\Q):= \widehat{\Gamma}(\Q) \cap \N^n$, called minimal generating set of $H(\Q)$, is a fundamental object in Weierstrass semigroup's Theory, see e.g. \cite{BQZ2018}, \cite{CT2005}, \cite{CT2018}, \cite{G2004}. More recently, the set  $\Lambda (\Q):= \widehat{\Lambda}(\Q) \cap \N^n$ has become a very important tool for determining the sets of gaps and pure gaps of $H (\Q)$, see \cite{CMT2024} and \cite{TT2019}. As a consequence of previous result, we can determine the sets $\Gamma(\Q)$ and $\Lambda(\Q)$ as well as its respective cardinalities as follows. We observe that the result below generalizes the results given in \cite[Propositions 4.6 and 4.7]{CMQ2024}. 

\begin{corollary}\label{coro_Upsilon}
Let $\Q=(Q_{1}, Q_{2}, \dots, Q_{n})$ and assume that $2\leq n \leq q$. Then
$$
	\Upsilon(\Q)=\Bigg\{(mj_1+t_{1}(i), \dots,  mj_n+t_{n}(i)):
	\begin{array}{l}
		1\leq i \leq m-1,\,\, j_1, \dots, j_n \in \N_0, \\
		j_1+\cdots +j_n=\beta(i)+1-n+\rho
	\end{array}
	\Bigg\},
	$$
	where 
$$\quad \rho=\left\{
\begin{array}{ll}
0, & \text{if }\, \Upsilon=\Gamma,\\
n-2, & \text{if }\, \Upsilon=\Lambda.
\end{array}\right.
$$  
In addition,
$$
|\Upsilon(\Q)|=\sum_{i=1}^{m-1}\binom{\be(i)+\rho}{n-1}.
$$
\end{corollary}
\begin{proof}
The description of $\Upsilon(\Q)$ follows directly from Theorem \ref{teo_widehat_Upsilon}, while the cardinality of $\Upsilon(\Q)$ is derived from the fact that the number of solutions to the equation $k_1 + \dots + k_n = k$ in non-negative integers is $\binom{k+n-1}{n-1}$, see \cite[Theorem 13.1]{LW2001}.
\end{proof}

In \cite{CMT2025}, a formula for the cardinality of the set $\Upsilon(\Q)$ is provided. To finish this section, we will show that this formula coincides with the formula given in Corollary \ref{coro_Upsilon}.

For $k$ a non-negative integer, define the set
$$
\Upsilon_{k, 0, \dots, 0}:=\Upsilon(\Q)\cap[km, (k+1)m)\times[0, m)^{n-1}.
$$ 
From the description of $\Upsilon(\Q)$ given in Corollary \ref{coro_Upsilon}, we obtain that 
$$
|\Upsilon_{k, 0, \dots, 0}|=|\{1\leq i\leq m-1: k=\be(i)+1-n+\rho\}|.
$$ 
Therefore,
\begin{align*}
|\Upsilon(\Q)|&=\sum_{i=1}^{m-1}\binom{\be(i)+\rho}{n-1}\\
&=\sum_{0\leq k}|\{1\leq i \leq m-1: \be(i)+\rho=k+n-1\}|\binom{k+n-1}{n-1}\\
&=\sum_{0\leq k}\binom{k+n-1}{n-1}|\Upsilon_{k, 0, \dots, 0}|.
\end{align*}
Thus, the formula for the cardinality of $\Upsilon(\Q)$ given in Corollary \ref{coro_Upsilon} coincides with the one provided in \cite[Proposition 3.2]{CMT2025}.


\section{Set of maximal elements for some function fields} \label{Section 4}

In this section, we apply our results to determine the sets $\widehat{\Upsilon}(\Q)$ and $\Upsilon(\Q)$ for places in certain function fields that have been studied recently in the theory of Weierstrass semigroups and coding theory. The first example is the function field of the maximal curves $\cX_{a, b, n, s}$ and $\cY_{n, s}$ introduced in \cite{TTT2016}, which in some cases cannot be covered by the Hermitian curve, and the second example is the function field of curves defined by $y^m=f(x)$, where $f(x) = (x - \al_1)\cdots (x - \al_t) \in \fq[x]$ and $\al_1, \al_2, \dots, \al_t \in \fq$ are pairwise distinct. In such function fields there is a characterization of the sets $\widehat{\Upsilon}(\Q)$ and $\Upsilon(\Q)$ with $\Q = (Q_{\infty}, Q_1, \ldots, Q_n)$, where $Q_{\infty}$ is the single place at infinity in the respective function field, see \cite{MT2023}, \cite{TT2019FF} and \cite{TT2019}. The techniques used in these works were based on the need for $Q_{\infty}$ belongs to the $n$-tuple $\Q$. For the second example, in \cite{CMT2025} and \cite{YH2017} the authors determined those sets for the case where $Q_{\infty}$ does not belong to the $n$-tuple $\Q$, with the condition $(t,m)=1$. Using the results presented in the previous section, we determine the sets $\widehat{\Upsilon}(\Q)$ and $\Upsilon(\Q)$ for the first example without requiring $Q_{\infty}$ to belong to the $n$-tuple $\Q$, and for the second example, without the condition $(t,m)=1$. In the third example, we determine the sets $\widehat{\Upsilon}(\Q)$ and $\Upsilon(\Q)$ for places in the function field of the Beelen-Montanucci maximal curve \cite{BM2018}. 

\subsection{Maximal curves that cannot be covered by the Hermitian curve}

Let $a, b, s \geq 1$ and $n \geq 3$ be integers such that $n$ is odd, $b$ divides $a$, and $b<a$. Let $q = p^a$ be a power of a prime $p$, let $c \in \fqs$ be such that $c^{q-1} = -1$, and suppose that $s$ divides $(q^n + 1)/(q + 1)$. Assuming these conditions, consider the following curves over $\fqsn$: 
$$
\cX_{a, b, n, s}:\quad cy^{\frac{q^n+1}{s}}=t(x)(t(x)^{q-1}+1)^{q+1}, \quad \text{where } t(x):=\textstyle\sum_{i=0}^{a/b-1}x^{p^{ib}},
$$
and 
$$
\cY_{n, s}: \quad y^{\frac{q^n+1}{s}}=(x^q+x)((x^q+x)^{q-1}-1)^{q+1}.
$$
These curves are maximal over $\fqsn$. From \cite[Theorem 3.5]{TTT2016}, the curve $\cX_{a,b,n,1}$ cannot be Galois-covered by the Hermitian curve $\cH_n$ over $\fqsn$. Furthermore, from \cite[Theorem 4.4]{TTT2016}, the curve $\cY_{3,s}$ cannot be covered by the Hermitian curve $\cH_3$ over $\fqss$ if $q >s(s+1)$.

Let $\cC$ be denote either the curve $\cX_{a, b, n, s}$ or the curve $\cY_{n, s}$. Let $Q_{\infty}$ be the only place at infinity in $\fqsn(\cC)$ and $Q_1, Q_2, \dots, Q_n$ be distinct pairwise places in $\fqsn(\cC)$ corresponding to roots of $t(x)$ if $\cC=\cX_{a, b, n, s}$, or to roots of $x^q+x$ if $\cC=\cY_{n, s}$. Moreover, let $\tilde{M}:=(q^n+1)/s$ and
$$
d:=\left\{
\begin{array}{ll}
	p^b, & \text{if } \cC=\cX_{a, b, n, s},\\
	1, & \text{if } \cC=\cY_{n, s}.
\end{array}\right.
$$

In \cite{MT2023}, the authors determine the sets $\widehat{\Upsilon}(\Q)$ and $\Upsilon(\Q)$, where $\Q=(Q_{\infty}, Q_1, \dots, Q_n)$ and $1\leq n \leq q/d$. In the following result, we obtain the sets $\widehat{\Upsilon}(\Q)$ and $\Upsilon(\Q)$ for $\Q=(Q_1, Q_2, \dots, Q_n)$ and $2\leq n \leq q/d$. 

\begin{proposition} 
Assume $2\leq n\leq q/d$ and let $\Q=(Q_1, \dots, Q_n)$. Then
\begin{align*}
	\widehat{\Upsilon}(\Q)=&\Bigg\{(\tilde{M} j_1+i, \dots,  \tilde{M} j_n+i):
	\begin{array}{l}
		1\leq i \leq \tilde{M} -1, \,\, j_1, \dots, j_n\in \Z, \\
		j_1+\cdots +j_n=\be(i)+1-n+\rho
	\end{array}
	\Bigg\}\\
	& \bigcup \Bigg\{\tilde{M} j_1, \dots, \tilde{M} j_n): j_1, \dots, j_n \in \Z, \, j_1+\cdots+j_n=\rho\Bigg\},
\end{align*}
where
$$
\be(i)=\frac{q}{d}+\frac{q(q-1)}{d}\ceil*{\frac{i(q+1)}{\tilde{M}}}-\floor*{\frac{iq^3}{d \tilde{M}}}-1\quad \text{and}\quad
\rho=\left\{
\begin{array}{ll}
	0, & \text{if } \Upsilon=\Gamma,\\
	n-2, & \text{if } \Upsilon=\Lambda.
\end{array}\right. 
$$
In particular,
$$
\Upsilon(\Q)=\Bigg\{(\tilde{M} j_1+i, \dots,  \tilde{M} j_n+i):
\begin{array}{l}
	 1\leq i \leq \tilde{M}-1, \,\,  j_1, \dots, j_n \in \N_0, \\
	 j_1+\cdots +j_n=\be(i)+1-n+\rho
\end{array}
\Bigg\}.
$$
\end{proposition}

\begin{proof}
First, note that the polynomials $t(x),\, t(x)^{q-1}+1,\, x^q+x$, and $(x^q + x)^{q-1} - 1$ are separable. Therefore the polynomials $t(x)(t(x)^{q-1}+1)^{q+1}$ and $(x^q+x)((x^q+x)^{q-1}-1)^{q+1}$ have $q/d$ zeros with multiplicity $1$, $q(q-1)/d$ zeros with multiplicity $q+1$, and one pole with multiplicity $q^3/d$ in $\fqsn(x)$.  
The result follows from the definition of $\be(i)$ given in (\ref{t beta}), Theorem \ref{teo_widehat_Upsilon}, and Corollary \ref{coro_Upsilon}. 
\end{proof}

\subsection{The curve $y^m=f(x)$, with $f(x)$ separable} Let $m, r \geq 2$ be integers such that $\char (\fq) \nmid m$, and let $f(x) = (x - \al_1)\cdots (x - \al_t) \in \fq[x]$, where $\al_1, \al_2, \dots, \al_t \in \fq$ are pairwise distinct. Consider the curve $\cX$ defined by the equation $y^m = f(x)$. For $1 \leq n \leq t$, let $Q_1, Q_2, \dots, Q_{n} \in \cP_{\fq(\cX)}$ denote the places corresponding to the zeros of $f(x)$. Moreover, if $(m, t) = 1$, let $Q_{\infty}$ denote the unique place at infinity in $\fq(\cX)$. We can find explicit descriptions for the sets $\widehat{\Gamma}(\Q)$, $\Gamma(\Q)$, $\widehat{\Lambda}(\Q)$, and $\Lambda(\Q)$, where $\Q=(Q_{\infty}, Q_1, \ldots , Q_n)$, in \cite[Theorem 6]{TT2019FF}, \cite[Theorem 9]{YH2017}, \cite[Corollary 4.3]{TT2019}, and \cite[Corollary 4.4]{TT2019}, respectively. For the case $\Q=(Q_1, \ldots , Q_n)$ with $(m,t)=1$ and $2\leq n \leq q$, we can find explicit descriptions of such sets in \cite[Theorem 5.6]{CMT2025}, \cite[Theorem 10]{YH2017}, \cite[Theorem 5.7]{CMT2025} and \cite[Corollary 5.8]{CMT2025}, respectively. All of these descriptions can be derived as particular cases of Theorem \ref{teo_widehat_Upsilon} and Corollary \ref{coro_Upsilon}.

Below we present the sets $\widehat{\Upsilon}(\Q)$ and $\Upsilon(\Q)$, where $\Q=(Q_1,\ldots,Q_n)$, without the condition $(m,t)=1$.

\begin{proposition} 
Assume $2\leq n \leq q$ and let $\Q=(Q_1, \dots, Q_n)$. Then
\begin{align*}
\widehat{\Upsilon}(\Q)=&\Bigg\{(mj_1+i, \dots,  mj_n+i):
\begin{array}{l}
	 1\leq i \leq m-1, \,\, j_1, \dots, j_n \in \Z, \text{ and}\\
	 j_1+\cdots +j_n=t-n-\floor*{it/m}+\rho
\end{array}
\Bigg\}\\
& \bigcup \Bigg\{(mj_1, \dots, mj_n): \, j_1, \dots, j_n \in \Z, \,\,  j_1+\cdots+j_n=\rho\Bigg\},
\end{align*}
where $\rho=0$ if $\Upsilon=\Gamma$ and $\rho=n-2$ if $\Upsilon=\Lambda$. In particular,
$$
\Upsilon(\Q)=\Bigg\{(mj_1+i, \dots, mj_n+i):
\begin{array}{l}
1\leq i \leq m-1, \,\, j_1, \dots, j_n\in \N_0, \text{ and}\\
j_1+\cdots +j_n=t-n-\floor*{it/m}+\rho
\end{array}
\Bigg\}.
$$
\end{proposition}
\begin{proof}
Note that the polynomial $f(x)$ has $t$ zeros with multiplicity $1$ and one pole with multiplicity $t$ in $\fq(x)$. 
Therefore, from the definition of $\be(i)$ given in (\ref{t beta}), we have $\beta(i)=t+\ceil*{-ti/m}-1=t-1-\floor*{ti/m}$. The result follows from Theorem \ref{teo_widehat_Upsilon} and Corollary \ref{coro_Upsilon}. 
\end{proof}

\subsection{Beelen-Montanucci curve}
The Beelen-Montanucci curve was introduced in \cite{BM2018} and is the second generalization of the $GK$ curve. A plane model for this curve was presented in \cite[Corollary 3.3]{MQ2022} as follows. For an odd integer $n \geq 3$, the Beelen-Montanucci curve over $\fqsn$ can be defined by the affine equation 
$$
\cB \cM _{n}:\quad y^{q^n+1}=(x^{q+1}-1)x^{q+1}(x^{q^2-1}-(x^{q+1}-1)^{q-1})^{q+1}.
$$
This curve is maximal over $\fqsn$, and for $n=3$ it is isomorphic to the $GK$ curve, which is the first example of a maximal curve not covered by the Hermitian curve.  In \cite{LV2022}, the authors determine the Weierstrass semigroup $H(P,Q)$ at certain two rational places in $\fqsn(\cB\cM_n)$.  

Let $Q_1, Q_2, \dots, Q_n$ be distinct pairwise places in $\fqsn(\cB\cM_n)$ corresponding to roots of the polynomial $x^{q+1}-1$ and let $M'=q^n+1$. In the following result we provided an explicit description of the sets $\widehat{\Upsilon}(\Q)$ and $\Upsilon(\Q)$ for $\Q=(Q_1, Q_2, \dots, Q_n)$. 

\begin{proposition} 
Assume $2 \leq n\leq q+1$ and let $\Q=(Q_1, \dots, Q_n)$. Then
\begin{align*}
\widehat{\Upsilon}(\Q)=&\Bigg\{(M'j_1+i, \dots,  M'j_n+i):
\begin{array}{l}
	1\leq i \leq M'-1, \,\, j_1, \dots, j_n\in \Z, \\
	j_1+\cdots +j_n=\be(i)+1-n+\rho
\end{array}
\Bigg\}\\
& \bigcup \Bigg\{(M'j_1, \dots, M'j_n): j_1, \dots, j_n \in \Z, \, j_1+\cdots+j_n=\rho\Bigg\},
\end{align*}
where
$$
\be(i)=q+1+(q^2-q-1)\ceil*{\frac{i(q+1)}{M'}}-\floor*{\frac{i(q^3-q)}{M'}}-1\quad \text{and}\quad
\rho=\left\{
\begin{array}{ll}
	0, & \text{if } \Upsilon=\Gamma,\\
	n-2, & \text{if } \Upsilon=\Lambda.
\end{array}\right. 
$$
In particular,
$$
\Upsilon(\Q)=\Bigg\{(M'j_1+i, \dots,  M'j_n+i):
\begin{array}{l}
	1\leq i \leq M'-1, \,\, j_1, \dots, j_n \in \N_0,\\
	j_1+\cdots +j_n=\be(i)+1-n+\rho
\end{array}
\Bigg\}.
$$
\end{proposition}

\begin{proof}
Let $g_1(x) = x^{q+1}-1$ and $g_2(x) = x^{q+1}(x^{q^2-1}-(x^{q+1}-1)^{q-1})^{q+1}$. Note that $g_1(x)$ has $q+1$ zeros with multiplicity $1$ and $g_2(x)$ has $q^2 - q -1$ zeros with multiplicity $q+1$. Moreover, $g(x):=g_1(x)g_2(x)$ has a unique pole with multiplicity $q^3-q$. The result follows from the definition of $\be(i)$ given in (\ref{t beta}), Theorem \ref{teo_widehat_Upsilon}, and Corollary \ref{coro_Upsilon}. 
\end{proof}


\bibliographystyle{abbrv}

\bibliography{generalizedws} 

\begin{thebibliography}{10}

\bibitem{ABQ2019}
M.~Abd\'{o}n, H.~Borges, and L.~Quoos.
\newblock Weierstrass points on {K}ummer extensions.
\newblock {\em Adv. Geom.}, 19(3):323--333, 2019.

\bibitem{ACGH1985}
E.~Arbarello, M.~Cornalba, P.~Griffiths, and J.~D. Harris.
\newblock {\em Geometry of Algebraic Curves}.
\newblock Springer-Verlag, 1985.

\bibitem{BMMQ2018}
D.~Bartoli, A.~M. Masuda, M.~Montanucci, and L.~Quoos.
\newblock Pure gaps on curves with many rational places.
\newblock {\em Finite Fields Appl.}, 53:287--308, 2018.

\bibitem{BQZ2018}
D.~Bartoli, L.~Quoos, and G.~Zini.
\newblock Algebraic geometric codes on many points from {K}ummer extensions.
\newblock {\em Finite Fields Appl.}, 52:319--335, 2018.

\bibitem{BMZ2020}
M.~Bartoli, D.~Montanucci and G.~Zini.
\newblock Weierstrass semigroups at every point of the suzuki curve.
\newblock {\em Acta Arithmetica}, 197:1--20, 2020.

\bibitem{BLM2021}
L.~Beelen, P.~Landi and M.~Montanucci.
\newblock Weierstrass semigroups on the skabelund maximal curve.
\newblock {\em Finite Fields Appl.}, 72:Paper No. 101811, 2021.

\bibitem{BM2018}
P.~Beelen and M.~Montanucci.
\newblock A new family of maximal curves.
\newblock {\em J. Lond. Math. Soc. (2)}, 98(3):573--592, 2018.

\bibitem{BMV2023}
P.~Beelen, M.~Montanucci, and L.~Vicino.
\newblock Weierstrass semigroups and automorphism group of a maximal curve with the third largest genus.
\newblock {\em Finite Fields Appl.}, 92:46--69, 2023.

\bibitem{BT2006}
P.~Beelen and N.~Tuta\c{s}.
\newblock A generalization of the {W}eierstrass semigroup.
\newblock {\em J. Pure Appl. Algebra}, 207(2):243--260, 2006.

\bibitem{CT2005}
C.~Carvalho and F.~Torres.
\newblock On {G}oppa codes and {W}eierstrass gaps at several points.
\newblock {\em Des. Codes Cryptogr.}, 35(2):211--225, 2005.

\bibitem{CT2018gk}
A.~Castellanos and G.~Tizziotti.
\newblock Weierstrass semigroup and pure gaps at several points on the gk curve.
\newblock {\em Bull Braz Math Soc}, 49:419--429, 2018.

\bibitem{CB2020}
A.~S. Castellanos and M.~Bras-Amor\'{o}s.
\newblock Weierstrass semigroup at {$m+1$} rational points in maximal curves which cannot be covered by the {H}ermitian curve.
\newblock {\em Des. Codes Cryptogr.}, 88(8):1595--1616, 2020.

\bibitem{CMQ2016}
A.~S. Castellanos, A.~M. Masuda, and L.~Quoos.
\newblock One- and two-point codes over {K}ummer extensions.
\newblock {\em IEEE Trans. Inform. Theory}, 62(9):4867--4872, 2016.

\bibitem{CMQ2024}
A.~S. Castellanos, E.~A.~R. Mendoza, and L.~Quoos.
\newblock Weierstrass semigroups, pure gaps and codes on function fields.
\newblock {\em Des. Codes Cryptography}, 92(5):1219--1242, 2024.

\bibitem{CMT2024}
A.~S. Castellanos, E.~A.~R. Mendoza, and G.~Tizziotti.
\newblock Complete set of pure gaps in function fields.
\newblock {\em J. Pure Appl. Algebra}, 228(4):Paper No. 107513, 20, 2024.

\bibitem{CMT2025}
A.~S. Castellanos, E.~A.~R. Mendoza, and G.~Tizziotti.
\newblock The set of pure gaps at several rational places in function fields.
\newblock {\em Des. Codes Cryptography}, pages 1--26, 2024.

\bibitem{CT2018}
A.~S. Castellanos and G.~Tizziotti.
\newblock On weierstrass semigroup at m points on curves of the form {f ( y )= g ( x )}.
\newblock {\em Journal of Pure and Applied Algebra}, 222(7):1803--1809, 2018.

\bibitem{C2008}
A.~D. Centina.
\newblock Weierstrass points and their impact in the study of algebraic curves: a historical account from the ``l\" uckensatz'' to the 1970s.
\newblock {\em Ann. Univ. Ferrara}, 54(1):37--59, 2008.

\bibitem{CMS2025}
E.~Cotterill, E.~A.~R. Mendoza, and P.~Speziali.
\newblock On gap sets in arbitrary {Kummer} extensions of {$K(x)$}.
\newblock In Preparation.

\bibitem{D1990}
F.~Delgado.
\newblock The symmetry of the {W}eierstrass generalized semigroups and affine embeddings.
\newblock {\em Proc. Amer. Math. Soc.}, 108(3):627--631, 1990.

\bibitem{DP2012}
I.~M. Duursma and S.~Park.
\newblock Delta sets for divisors supported in two points.
\newblock {\em Finite Fields Appl.}, 18(5):865--885, 2012.

\bibitem{GKL1993}
A.~Garcia, K.~S. Jeong, and R.~F. Lax.
\newblock Consecutive {Weierstrass} gaps and minimum distance of {Goppa} codes.
\newblock {\em J. Pure Appl. Algebra}, 84(2):199--207, 1993.

\bibitem{HLP1998}
T.~Hoholdt, J.~van Lint, and R.~Pellikaan.
\newblock {\em Algebraic geometry codes}.
\newblock Handbook of Coding Theory, vol. 1. V.S. Pless and W.C. Huffman, first edition, 1998.

\bibitem{H1996}
M.~Homma.
\newblock The {W}eierstrass semigroup of a pair of points on a curve.
\newblock {\em Arch. Math. (Basel)}, 67(4):337--348, 1996.

\bibitem{K1994}
S.~J. Kim.
\newblock On the index of the {W}eierstrass semigroup of a pair of points on a curve.
\newblock {\em Arch. Math. (Basel)}, 62(1):73--82, 1994.

\bibitem{KNT2020}
G.~Korchm{\'a}ros, G.~P. Nagy, and M.~Timpanella.
\newblock Codes and gap sequences of {Hermitian} curves.
\newblock {\em IEEE Trans. Inf. Theory}, 66(6):3547--3554, 2020.

\bibitem{LV2022}
L.~Landi and L.~Vicino.
\newblock Two-point {AG} codes from the {B}eelen-{M}ontanucci maximal curve.
\newblock {\em Finite Fields Appl.}, 80:Paper No. 102009, 17, 2022.

\bibitem{M2004}
H.~Maharaj.
\newblock Code construction on fiber products of {K}ummer covers.
\newblock {\em IEEE Trans. Inform. Theory}, 50(9):2169--2173, 2004.

\bibitem{GD2010}
G.~Mattews and J.~D. Peachey.
\newblock Minimal generating sets of weierstrass semigroups of certain $m$-tuples on the norm-trace function field.
\newblock {\em Contemporary Mathematics}, 518:315--326, 2010.

\bibitem{GM2001}
G.~L. Matthews.
\newblock Weierstrass pairs and minimum distance of {G}oppa codes.
\newblock {\em Des. Codes Cryptogr.}, 22(2):107--121, 2001.

\bibitem{GM2004}
G.~L. Matthews.
\newblock Codes from the suzuki function field.
\newblock {\em IEEE Trans. Inform. Theory}, 50(12):3298--3302, 2004.

\bibitem{G2004}
G.~L. Matthews.
\newblock The {W}eierstrass semigroup of an {$m$}-tuple of collinear points on a {H}ermitian curve.
\newblock In {\em Finite fields and applications}, volume 2948 of {\em Lecture Notes in Comput. Sci.}, pages 12--24. Springer, Berlin, 2004.

\bibitem{E2023}
E.~A.~R. Mendoza.
\newblock On {Kummer} extensions with one place at infinity.
\newblock {\em Finite Fields Appl.}, 89:102209, 2023.
\newblock Id/No 102209.

\bibitem{MQ2022}
E.~A.~R. Mendoza and L.~Quoos.
\newblock Explicit equations for maximal curves as subcovers of the {{\(BM\)}} curve.
\newblock {\em Finite Fields Appl.}, 77:22, 2022.
\newblock Id/No 101945.

\bibitem{MT2023}
M.~Montanucci and G.~Tizziotti.
\newblock Generalized {Weierstrass} semigroups at several points on certain maximal curves which cannot be covered by the {Hermitian} curve.
\newblock {\em Des. Codes Cryptography}, 91(3):831--851, 2023.

\bibitem{MTT2019}
J.~J. Moyano-Fern{\'a}ndez, W.~Ten{\'o}rio, and F.~Torres.
\newblock Generalized {Weierstrass} semigroups and their {Poincar{\'e}} series.
\newblock {\em Finite Fields Appl.}, 58:46--69, 2019.

\bibitem{MTT2008}
C.~Munuera, G.~Tizziotti, and F.~Torres.
\newblock Two-points codes on norm-trace curves.
\newblock {\em Second International Castle Meeting, ICMCTA 2008 (A. Barbero Ed.), Lecture Notes in Comput. Sci., Springer-Verlag Berlin Heidelberg}, 5228:128--136, 2008.

\bibitem{S2009}
H.~Stichtenoth.
\newblock {\em Algebraic function fields and codes}, volume 254 of {\em Graduate Texts in Mathematics}.
\newblock Springer-Verlag, Berlin, second edition, 2009.

\bibitem{TTT2016}
S.~Tafazolian, A.~Teher{\'a}n-Herrera, and F.~Torres.
\newblock Further examples of maximal curves which cannot be covered by the {Hermitian} curve.
\newblock {\em J. Pure Appl. Algebra}, 220(3):1122--1132, 2016.

\bibitem{TT2019FF}
W.~Ten\'{o}rio and G.~Tizziotti.
\newblock Generalized {W}eierstrass semigroups and {R}iemann-{R}och spaces for certain curves with separated variables.
\newblock {\em Finite Fields Appl.}, 57:230--248, 2019.

\bibitem{TT2019}
W.~Ten\'{o}rio and G.~Tizziotti.
\newblock On {W}eierstrass gaps at several points.
\newblock {\em Bull. Braz. Math. Soc. (N.S.)}, 50(2):543--559, 2019.

\bibitem{LW2001}
J.~van Lint and R.~Wilson.
\newblock {\em A course in combinatorics}.
\newblock Cambridge University Press, Cambridge, second edition, 2001.

\bibitem{YH2017}
S.~Yang and C.~Hu.
\newblock Weierstrass semigroups from {K}ummer extensions.
\newblock {\em Finite Fields Appl.}, 45:264--284, 2017.

\end{thebibliography}

\end{document}